\documentclass{article}
\usepackage{amsfonts}
\usepackage{latexsym}
\usepackage{mathrsfs}
\usepackage{amssymb}
\usepackage{amsmath}
\usepackage{amsthm}
\usepackage{indentfirst}
\usepackage{color}
\usepackage{float}
\usepackage{graphicx}

\allowdisplaybreaks
\newtheorem{theorem}{\color{black}\indent Theorem}[section]

\newtheorem{definition}{\color{black}\indent Definition}[section]
\newtheorem{remark}{\color{black}\indent Remark}[section]

\begin{document}
\title{\LARGE\bf Global asymptotic behavior of solutions to a class of Kirchhoff equations}
\author{Yuzhu Han$^\dag$}
 \date{}
 \maketitle

 \footnotetext{\hspace{-1.9mm}$^\dag$Corresponding author.\\
  Email addresses: yzhan@jlu.edu.cn(Y. Han).

\thanks{
$^*$The project is supported by NSFC (11401252),
by Science and Technology Development Project of Jilin Province (20160520103JH) 
and by The Education Department of Jilin Province (JJKH20190018KJ).}}

\begin{center}
{\it\small School of Mathematics, Jilin University, Changchun 130012, P.R. China}
\end{center}

\date{}
\maketitle

{\bf Abstract}\ In this paper, a parabolic type Kirchhoff equation and its stationary counterpart are considered.
For the evolution problem, the precise decay rates of the weak solution and of the corresponding energy functional are derived.
For the stationary problem, a ground-state solution is obtained by applying Lagrange multiplier method.
Moreover, the asymptotic behaviors of the general global solutions are also described.
These results extend some recent ones obtained in [Threshold results for the existence of global and blow-up solutions to
Kirchhoff equations with arbitrary initial energy, Computers and Mathematics with Applications, 75(2018), 3283-3297] by Han and Li.



\section{Introduction}
\setcounter{equation}{0}

In this paper, we are interested in the following initial boundary value problem (IBVP) for a class of nonlocal parabolic
equation
\begin{equation}\label{1.1}
\begin{cases}
u_t-M(\int_{\Omega}|\nabla u|^2\mathrm{d}x)\Delta u=|u|^{q-1}u, &(x,t)\in \Omega\times(0,T),\\
u=0, &(x,t)\in \partial \Omega\times(0,T),\\
u(x,0)=u_0(x), & x\in\Omega,
\end{cases}
\end{equation}
where $M(s)=a+bs$ with $a, b>0$, $\Omega\subset\mathbb{R}^n(n\geq1)$ is a bounded smooth domain with $\partial\Omega$ as its boundary,
$3<q<2^*-1$, where $2^*$ is the Sobolev conjugate of $2$, i.e. $2^*=+\infty$ for $n=1,2$ and $2^*=\dfrac{2n}{n-2}$ for $n\geq 3$.
Moreover, the initial datum $u_0$ belongs to the energy space $H_0^1(\Omega)$.

When he describes the transversal oscillations of a stretched string by taking into account the change in the string length,
Kirchhoff \cite{kirchhoff} first proposed the hyperbolic counterpart of the equation in \eqref{1.1}.
The main feature of such equations is that the coefficient of the second order term depends on
the $L^2(\Omega)$ norm of the gradient of the unknown, which means they are no longer pointwise identities.
Therefore, such equations are usually referred to as Kirchhoff equations or nonlocal equations.
When the nonlinearity $|u|^{q-1}u$ in \eqref{1.1} is replaced by some $L^2(\Omega)$ function $f(x)$,
Chipot et al. \cite{Chipot2003} studied the local and global existence, uniqueness and global asymptotic behavior of weak solutions to \eqref{1.1}.
Recently, Han et al. \cite{Han-Li,Han-Gao} considered the global existence and finite time blow-up of solutions to IBVP \eqref{1.1}
by applying the modified potential well method, which was first proposed by Sattinger \cite{Sattinger1968}, Payne and Sattinger \cite{Payne1975},
improved by Liu \cite{YCLiu2003} and Xu et al. \cite{RZXu2013} and then applied to other evolution problems, for example, in \cite{Han,Qu}.

The purpose of this paper is to generalize some results obtained in \cite{Han-Li}.
As in \cite{Han-Li}, we will denote by $\|u\|_{r}$ the $L^r(\Omega)$ norm of a Lebesgue function $u\in L^r(\Omega)$ for $r\geq 1$ and by $(\cdot, \cdot)$ the inner product in $L^2(\Omega)$.
By $H_0^1(\Omega)$ we denote the Sobolev space such that both $u$ and $|\nabla u|$ belong to $L^2(\Omega)$ for any $u\in H_0^1(\Omega)$.
Recalling Poincar\'{e}'s inequality, we can equip $H_0^1(\Omega)$ with the equivalent norm $\|u\|_{H_0^1(\Omega)}=\|\nabla u\|_2$.
We denote by $H^{-1}(\Omega)$ the dual space to $H_0^1(\Omega)$.
For $u\in H_0^1(\Omega)$, set
\begin{equation}\label{1.2}
J(u)=\dfrac{a}{2}\|\nabla u\|_2^2+\dfrac{b}{4}\|\nabla u\|_2^4-\dfrac{1}{q+1}\|u\|^{q+1}_{q+1},
\end{equation}
\begin{equation}\label{1.3}
I(u)=a\|\nabla u\|_2^2+b\|\nabla u\|_2^4-\|u\|^{q+1}_{q+1},
\end{equation}
and define the Nehari manifold
\begin{equation}\label{1.4}
\mathcal{N}=\{u\in H_0^1(\Omega)| \ I(u)=0, \ \|\nabla u\|_2\neq0\}.
\end{equation}
Since $q+1<2^*$, it can be directly checked that both $J(u)$ and $I(u)$ are well defined in $H_0^1(\Omega)$
and $J, I\in C^1(H_0^1(\Omega), \mathbb{R})$. Moreover, for any $u,v\in H_0^1(\Omega)$,
\begin{equation}\label{1.5}
\langle J'(u),v\rangle=a(\nabla u,\nabla v)+b\|\nabla u\|_2^2(\nabla u,\nabla v)-(|u|^{q-1}u,v),
\end{equation}
and
\begin{equation}\label{1.6}
\langle I'(u),v\rangle=2a(\nabla u,\nabla v)+4b\|\nabla u\|_2^2(\nabla u,\nabla v)-(q+1)(|u|^{q-1}u,v),
\end{equation}
where $\langle\cdot,\cdot\rangle$ denotes the pairing between $H^{-1}(\Omega)$ and $H_0^1(\Omega)$.
From \eqref{1.2} and \eqref{1.3} it follows that
\begin{equation}\label{1.7}
\begin{split}
J(u)&=\dfrac{a(q-1)}{2(q+1)}\|\nabla u\|_2^2+\dfrac{b(q-3)}{4(q+1)}\|\nabla u\|_2^4+\dfrac{1}{q+1}I(u)\\
&=\dfrac{a}{4}\|\nabla u\|_2^2+\dfrac{q-3}{4(q+1)}\|u\|_{q+1}^{q+1}+\dfrac{1}{4}I(u).
\end{split}
\end{equation}
The potential well and its corresponding set are defined, respectively, by
\begin{eqnarray*}
&&W=\{u\in H_0^1(\Omega)| \ I(u)>0, \ J(u)<d\}\cup\{0\},\\
&&V=\{u\in H_0^1(\Omega)| \ I(u)<0, \ J(u)<d\},
\end{eqnarray*}
where
\begin{equation}\label{d}
d=\inf_{0\neq u\in H_0^1(\Omega)}\sup_{\lambda>0}J(\lambda u)=\inf_{u\in \mathcal{N}}J(u)
\end{equation}
is the depth of the potential well $W$. By Lemma 2.1 in \cite{Han-Li},
\begin{equation}\label{d0}
d\geq d_0:=\dfrac{a(q-1)}{2(q+1)}(\dfrac{a}{S^{q+1}})^{\frac{2}{q-1}},
\end{equation}
where $S>0$ is the optimal embedding constant from $H_0^1(\Omega)$ to $L^{q+1}(\Omega)$,
i.e.,
\begin{equation}\label{s}
\dfrac{1}{S}=\inf\limits_{0\neq w\in H_0^1(\Omega)}\dfrac{\|\nabla w\|_2}{\|w\|_{q+1}}.
\end{equation}

In view of the definition of $d$, a natural question is that can $d$ be attained?
We shall give a positive answer to this question and show that it is closely related to the stationary problem of IBVP \eqref{1.1},
i.e., the following boundary value problem (BVP)
\begin{equation}\label{bvp}
\begin{cases}
-M(\int_{\Omega}|\nabla u|^2\mathrm{d}x)\Delta u=|u|^{q-1}u, &x\in \Omega,\\
u=0, &x\in \partial \Omega.
\end{cases}
\end{equation}
It is well known that $u\in H_0^1(\Omega)$ is a weak solution to \eqref{bvp} if and only if it is a critical point of $J$.
Let $\mathcal{S}$ denote the set of all weak solutions to BVP \eqref{bvp}.
Then we have the following theorem.

\begin{theorem}\label{stationary}
Let $\mathcal{N}$ and $d$ be defined in \eqref{1.4} and \eqref{d}, respectively.
Then there is a $v_0\in \mathcal{N}$ such that

(1) $J(v_0)=\inf\limits_{v\in\mathcal{N}}J(v)=d;$

(2) $v_0$ is a ground-state solution to BVP \eqref{bvp}, i.e., $v_0\in \mathcal{S}\setminus\{0\}$ and
$$J(v_0)=\inf\limits_{v\in\mathcal{S}\setminus\{0\}}J(v).$$
\end{theorem}

\begin{remark}
Although we stated Theorem \ref{stationary} when $3<q<2^*-1$, it is also valid when $q=3$,
provided we assume in addition that $0<b<S^4$, where $S>0$ is given in \eqref{s}.
In fact, noticing that the constant $S>0$ in \eqref{s} is attained,
it is easily checked that the Nehari's manifold $\mathcal{N}$ is non-empty and the potential depth $d$ is still positive in this case.
The remaining proof is similar to that of Theorem \ref{stationary}, and interested readers may check the details.
\end{remark}

The weak solution to IBVP \eqref{1.1} is defined as follows.
\begin{definition}\label{de2.1}$\mathrm{\bf{(Weak \ solution)}}$
Let $T>0$. A function $u=u(x,t)\in L^{\infty}(0, T; H_0^1(\Omega))$ with $u_t\in L^2(0, T; L^2(\Omega))$ is called a weak solution to IBVP \eqref{1.1} on $\Omega\times[0,T)$, if
$u(x, 0)=u_0\in H_0^1(\Omega)$ and satisfies
\begin{equation}\label{weak-solution}
(u_t, \phi)+\Big((a+b\int_{\Omega}|\nabla u|^2\mathrm{d}x)\nabla u, \nabla\phi\Big)=(|u|^{q-1}u, \phi), \ \ \ \ a.\ e. \ t\in(0, T),\\
\end{equation}
for any $\phi\in H_0^1(\Omega)$. Moreover, $u(x,t)$ satisfies
\begin{equation}\label{ei}
\int_0^t\|u_{\tau}\|_2^2\mathrm{d}\tau+J(u(x,t))=J(u_0), \qquad \ a.\ e. \ t\in(0, T).
\end{equation}
\end{definition}

Sometimes, we think of the function $u(x,t)$ of the space variable $x$ as an element of $H_0^{1}(\Omega)$,
and briefly denote the element of $H_0^{1}(\Omega)$ that arises this way by $u(t)$; therefore $u(t)\in H_0^{1}(\Omega)$.
When $J(u_0)<d$, the existence of global solutions to IBVP \eqref{1.1}, which can be summarized in the following theorem, was obtained in \cite{Han-Li}.

\begin{theorem}\label{low} [3, Theorem 3.1]
Assume $a, b >0$, $3<q<2^*-1$ and $u_0\in H_0^1(\Omega)$. If
$J(u_0)<d$ and $I(u_0)>0$, then IBVP \eqref{1.1} admits a global weak solution $u\in L^{\infty}(0,\infty; H_0^1(\Omega))$
with $u_t\in L^2(0,\infty; L^2(\Omega))$ and $u(t)\in W$ for $0\leq t<\infty$.
Moreover, $\|u\|_2^2\leq\|u_0\|_2^2e^{-Ct}$,
where $C>0$ is a constant. In addition, the weak solution is unique if it is bounded.
\end{theorem}

In view of Theorem \ref{low}, we observe that the following three problems are unsolved:

(i) Can we give the exact value or can we find some positive bounds for the exponential rate $C$ in Theorem \ref{low}?

(ii) Only the decay rate of $\|u(t)\|_2^2$ was derived. Can we also estimate the decay rate of $\|\nabla u(t)\|_2^2$ and $\|u(t)\|_{q+1}^2$?

(iii) It follows from \eqref{ei} that $J(u(t))$ is decreasing with respect to $t$. Can we derive the decay rate of $J(u(t))$?

The second result in this paper is concerned with the above problems, which can be solved for the case of $J(u_0)<d_0\leq d$.

\begin{theorem}\label{rate}
Assume $a, b >0$, $3<q<2^*-1$ and $u_0\in H_0^1(\Omega)$. If
$J(u_0)<d_0$ and $I(u_0)>0$, then IBVP \eqref{1.1} admits a global weak solution $u\in L^{\infty}(0,\infty; H_0^1(\Omega))$
with $u_t\in L^2(0,\infty; L^2(\Omega))$ and $u(t)\in W$ for $0\leq t<\infty$.
Moreover,
$$\|u\|_2^2\leq\|u_0\|_2^2e^{-C_1t},$$
$$\|\nabla u(t)\|_2^2\leq \dfrac{2(q+1)}{a(q-1)}[J(u_0)+\|u_0\|_2^2]e^{-C_2t},$$
$$\|u(t)\|_{q+1}^2\leq \dfrac{2S^2(q+1)}{a(q-1)}[J(u_0)+\|u_0\|_2^2]e^{-C_2t},$$
$$J(u(t))+\|u(t)\|_2^2\leq [J(u_0)+\|u_0\|_2^2]e^{-C_2t},$$
where $C_1=2a\lambda_1\Big[1-\big(\dfrac{J(u_0)}{d_0}\big)^{\frac{q-1}{2}}\Big]>0$,
$C_2=\dfrac{a\lambda_1\alpha(q-1)}{a\lambda_1(q-1)+2(q+1)}>0$,
 $\lambda_1>0$ is the first eigenvalue of $-\Delta$ in $\Omega$ under homogeneous Dirichlet boundary condition and $\alpha=8\Big[1-\big(\dfrac{J(u_0)}{d_0}\big)^{\frac{q-1}{2}}\Big]\Big[2-\dfrac{4}{q+1}\big(\dfrac{J(u_0)}{d_0}\big)^{\frac{q-1}{2}}\Big]^{-1}>0$.
\end{theorem}

\begin{remark}
From \eqref{1.7} and $I(u_0)>0$ one sees that $J(u_0)>0$. Therefore, all the terms in Theorem \ref{rate} make sense.
\end{remark}

In \cite{Han-Li}, the authors also investigated the global existence and finite time blow-up of solutions to IBVP \eqref{1.1} when $J(u_0)>d$.
To introduce these results, set
\begin{eqnarray*}
&& \mathcal{N}_+=\{u\in H_0^1(\Omega)| \ I(u)>0\},\quad \mathcal{N}_-=\{u\in H_0^1(\Omega)| \ I(u)<0\},\\
&&  J^s=\{u\in H_0^1(\Omega)| \ J(u)<s\},\quad \mathcal{N}_s=\mathcal{N}\cap J^s.
\end{eqnarray*}
Then for any $s>d$,
\begin{eqnarray}\label{52.1}
&& \mathcal{N}_s=\Big\{u\in \mathcal{N}:\dfrac{a(q-1)}{2(q+1)}\|\nabla u\|_2^2+\frac{b(q-3)}{4(q+1)}\|\nabla u\|_2^4<s\Big\}\neq\emptyset.
\end{eqnarray}
Define
\begin{equation}\label{52.2}
\lambda_s:=\inf\{\|u\|_2\mid u\in \mathcal{N}_s\},\quad\Lambda_s:=\sup\{\|u\|_2\mid u\in \mathcal{N}_s\}.
\end{equation}

With the help of the notations given above, Han and Li gave another global existence and finite time blow-up result
for $J(u_0)>d$.

\begin{theorem}\label{high} [3, Theorem 5.1]
Let $3<q<2^*-1$. Assume that $J(u_0)>d$, then the following statements hold.

(i) If $u_0\in \mathcal{N}_+$ and $\|u_0\|_2\leq\lambda_{J(u_0)}$,
then all the solutions to IBVP \eqref{1.1} exist globally and tend to zero in $H_0^1(\Omega)$ as $t\rightarrow\infty$;

(ii) If $u_0\in \mathcal{N}_-$ and $\|u_0\|_2\geq\Lambda_{J(u_0)}$,
then all the solutions to IBVP \eqref{1.1} blow up in finite time.
\end{theorem}

To make Theorem \ref{high} nontrivial, it is necessary to show that $K_1\leq \lambda_{J(u_0)}\leq\Lambda_{J(u_0)}\leq K_2$
for some positive constants $K_1$ and $K_2$, which was not done in \cite{Han-Li}.
Moreover, from Theorem \ref{high} we know that the global solutions converge to $0$ in $H_0^1(\Omega)$
as $t\rightarrow\infty$ when the initial data satisfy some specific conditions.
Can we describe the asymptotic behaviors of the general global solutions?
At the end of this section, we shall answer these two questions.
Our results in these two directions can be summarized into the following two theorems.

\begin{theorem}\label{bound}
Let $3<q<2^*-1$. Then for any $s>d$, $\lambda_s$ and $\Lambda_s$ defined in \eqref{52.2} satisfy
\begin{equation*}
K_1\leq\lambda_s\leq\Lambda_s\leq K_2,
\end{equation*}
where
\begin{equation*}
K_1=
\begin{cases}
\Big(\dfrac{a}{G}\Big)^{1/\gamma}\theta^{\frac{4-n(q-1)}{2\gamma}},\qquad if\ q\leq 1+\dfrac{4}{n};\\
\Big(\dfrac{a}{G}\Big)^{1/\gamma}\Big[\dfrac{2(q+1)s}{a(q-1)}\Big]^{\frac{4-n(q-1)}{4\gamma}},\qquad if\ q>1+\dfrac{4}{n},
\end{cases}
\end{equation*}
$\gamma=q+1-\dfrac{n(q-1)}{2}>0$, $\theta= \Big(\dfrac{2(q+1)d}{q-1}\Big)^{\frac{1}{q+1}}S^{-1}$,
$G>0$ is the constant in Gagliardo-Nirenberg inequality and $K_2=\sqrt{\dfrac{2(q+1)s}{a\lambda_1(q-1)}}$.
\end{theorem}

\begin{theorem}\label{asy}
Assume $a, b >0$, $3<q<2^*-1$ and $u_0\in H_0^1(\Omega)$. Let $u=u(t)$
be a global solution to IBVP \eqref{1.1}. Then there exists a $u^*\in \mathcal{S}$
and an increasing sequence $\{t_k\}_{k=1}^\infty$ with $t_k\rightarrow\infty$ as $k\rightarrow\infty$ such that
$$\lim\limits_{k\rightarrow\infty}\|u(t_k)-u^*\|_{H_0^1(\Omega)}=\lim\limits_{k\rightarrow\infty}\|\nabla u(t_k)-\nabla u^*\|_2=0.$$
\end{theorem}

\par
\section{Proofs of the main results.}
\setcounter{equation}{0}

\textbf{Proof of Theorem \ref{stationary}.} \emph{(1)} From \eqref{1.7} and the definitions of $d$ and $\mathcal{N}$ it follows that
\begin{equation*}
d=\inf_{v\in \mathcal{N}}J(v)=\inf_{v\in \mathcal{N}}\Big[\dfrac{a(q-1)}{2(q+1)}\|\nabla v\|_2^2+\dfrac{b(q-3)}{4(q+1)}\|\nabla v\|_2^4\Big].
\end{equation*}
Let $\{v_k\}_{k=1}^\infty\subset\mathcal{N}$ be a minimizing sequence of $J$. Then
\begin{equation}\label{2.1}
\lim\limits_{k\rightarrow\infty}J(v_k)=\lim\limits_{k\rightarrow\infty}\Big[\dfrac{a(q-1)}{2(q+1)}\|\nabla v_k\|_2^2+\dfrac{b(q-3)}{4(q+1)}\|\nabla v_k\|_2^4\Big]=d.
\end{equation}
Noticing that $q>3$, we obtain from \eqref{2.1} that $\{v_k\}_{k=1}^\infty\subset\mathcal{N}$ is bounded in $H_0^1(\Omega)$,
which, together with the fact that $q+1<2^*$, implies that there is a subsequence of $\{v_k\}_{k=1}^\infty\subset\mathcal{N}$,
which we still denote by $\{v_k\}_{k=1}^\infty\subset\mathcal{N}$, and a $v_0\in H_0^1(\Omega)$ such that
\begin{equation}\label{2.2}
\begin{split}
v_k\rightharpoonup v_0 \ & weakly \ in \ H_0^1(\Omega) \ as\ k \rightarrow\infty,\\
v_k\rightarrow v_0 \ & \ strongly \ in \ L^{q+1}(\Omega) \ as\ k \rightarrow\infty.
\end{split}
\end{equation}
Since $\{v_k\}_{k=1}^\infty\subset\mathcal{N}$, we have
$$a\|\nabla v_k\|_2^2+b\|\nabla v_k\|_2^4=\|v_k\|^{q+1}_{q+1}.$$
Combining this identity with the weakly lower semi-continuity of $\|\cdot\|_2$ and recalling \eqref{2.2} we get
\begin{equation}\label{2.3}
\begin{split}
a\|\nabla v_0\|_2^2+b\|\nabla v_0\|_2^4&\leq\liminf_{k\rightarrow\infty}[a\|\nabla v_k\|_2^2+b\|\nabla v_k\|_2^4]\\
&=\lim_{k\rightarrow\infty}[a\|\nabla v_k\|_2^2+b\|\nabla v_k\|_2^4]\\
&=\lim_{k\rightarrow\infty}\|v_k\|^{q+1}_{q+1}=\|v_0\|^{q+1}_{q+1}.
\end{split}
\end{equation}

We claim that $I(v_0)=0$, or equivalently $\|\nabla v_0\|_2=\lim\limits_{k\rightarrow\infty}\|\nabla v_k\|_2$.
If not, then by \eqref{2.3} it must hold that $a\|\nabla v_0\|_2^2+b\|\nabla v_0\|_2^4<\|v_0\|^{q+1}_{q+1}$ and $\|\nabla v_0\|_2\neq0$.
By Lemma 2.2 in \cite{Han-Li}, there exists a unique $\lambda^*>0$ such that $\lambda^* v_0\in\mathcal{N}$, i.e., $I(\lambda^* v_0)=0$.
Therefore,
\begin{equation*}
a\lambda^{*2}\|\nabla v_0\|_2^2+b\lambda^{*4}\|\nabla v_0\|_2^4=\lambda^{*q+1}\|v_0\|^{q+1}_{q+1}>\lambda^{*q+1}(a\|\nabla v_0\|_2^2+b\|\nabla v_0\|_2^4),
\end{equation*}
i.e.,
\begin{equation*}
a(\lambda^{*q+1}-\lambda^{*2})\|\nabla v_0\|_2^2+b(\lambda^{*q+1}-\lambda^{*4})\|\nabla v_0\|_2^4<0.
\end{equation*}
This implies that $\lambda^*<1$.

On one hand, we combine \eqref{2.1} with the fact that $\|\nabla v_0\|_2<\lim\limits_{k\rightarrow\infty}\|\nabla v_k\|_2$ to obtain
\begin{equation*}
\begin{split}
&\dfrac{a(q-1)}{2(q+1)}\|\nabla v_0\|_2^2+\dfrac{b(q-3)}{4(q+1)}\|\nabla v_0\|_2^4\\
<&\dfrac{a(q-1)}{2(q+1)}\lim\limits_{k\rightarrow\infty}\|\nabla v_k\|_2^2+\dfrac{b(q-3)}{4(q+1)}\lim\limits_{k\rightarrow\infty}\|\nabla v_k\|_2^4\\
=&\lim\limits_{k\rightarrow\infty}\Big[\dfrac{a(q-1)}{2(q+1)}\|\nabla v_k\|_2^2+\dfrac{b(q-3)}{4(q+1)}\|\nabla v_k\|_2^4\Big]\\
=&d.
\end{split}
\end{equation*}

On the other hand, it follows from $I(\lambda^* v_0)=0$, $\lambda^*<1$ and \eqref{1.7} that
\begin{equation*}
\begin{split}
d\leq& J(\lambda^* v_0)=\dfrac{a(q-1)}{2(q+1)}\lambda^{*2}\|\nabla v_0\|_2^2+\dfrac{b(q-3)}{4(q+1)}\lambda^{*4}\|\nabla v_0\|_2^4\\
<&\dfrac{a(q-1)}{2(q+1)}\|\nabla v_0\|_2^2+\dfrac{b(q-3)}{4(q+1)}\|\nabla v_0\|_2^4,
\end{split}
\end{equation*}
a contradiction. Therefore, $I(v_0)=0$ and $\|\nabla v_0\|_2=\lim\limits_{k\rightarrow\infty}\|\nabla v_k\|_2$.
Together with the weak convergence in \eqref{2.2} and the uniform convexity of $H_0^1(\Omega)$, it implies that
\begin{equation*}
v_k\rightarrow v_0 \ \ strongly \ in \ H_0^1(\Omega) \ as\ k \rightarrow\infty.
\end{equation*}
Moreover,
\begin{equation*}
\begin{split}
J(v_0)&=\dfrac{a(q-1)}{2(q+1)}\|\nabla v_0\|_2^2+\dfrac{b(q-3)}{4(q+1)}\|\nabla v_0\|_2^4+\dfrac{1}{q+1}I(v_0)\\
&=\dfrac{a(q-1)}{2(q+1)}\|\nabla v_0\|_2^2+\dfrac{b(q-3)}{4(q+1)}\|\nabla v_0\|_2^4\\
&=\lim\limits_{k\rightarrow\infty}\Big[\dfrac{a(q-1)}{2(q+1)}\|\nabla v_k\|_2^2+\dfrac{b(q-3)}{4(q+1)}\|\nabla v_k\|_2^4\Big]\\
&=d,
\end{split}
\end{equation*}
which implies that $v_0\neq0$. Therefore $v_0\in\mathcal{N}$ and $J(v_0)=\inf\limits_{v\in\mathcal{N}}J(v)=d$.

\emph{(2)} Since $J(v_0)=\inf\limits_{v\in\mathcal{N}}J(v)=d$, by the theory of Lagrange multipliers,
there exists a constant $\sigma\in\mathbb{R}$ such that
\begin{equation}\label{2.4}
J'(v_0)-\sigma I'(v_0)=0,
\end{equation}
which then ensures
\begin{equation}\label{2.5}
\sigma \langle I'(v_0),v_0\rangle=\langle J'(v_0),v_0\rangle=I(v_0)=0.
\end{equation}
Recalling \eqref{1.6} and the fact that $I(v_0)=0$ we obtain
\begin{equation*}
\begin{split}
\langle I'(v_0),v_0\rangle&=2a\|\nabla v_0\|_2^2+4b\|\nabla v_0\|_2^4-(q+1)\|v_0\|_{q+1}^{q+1}\\
&=2a\|\nabla v_0\|_2^2+4b\|\nabla v_0\|_2^4-(q+1)[a\|\nabla v_0\|_2^2+b\|\nabla v_0\|_2^4]\\
&=-a(q-1)\|\nabla v_0\|_2^2-b(q-3)\|\nabla v_0\|_2^4<0.
\end{split}
\end{equation*}
Therefore, $\sigma=0$, which in turn implies $J'(v_0)=0$ by \eqref{2.4}. So $v_0\in \mathcal{S}\setminus\{0\}$.
Noticing that $\mathcal{S}\setminus\{0\}\subset\mathcal{N}$, we have
$J(v_0)=\inf\limits_{v\in\mathcal{S}\setminus\{0\}}J(v)$ and the proof of Theorem \ref{stationary} is complete. \quad\qquad$\Box$

\textbf{Proof of Theorem \ref{rate}.} Since $d_0\leq d$,
it follows from Theorem \ref{low} that IBVP \eqref{1.1} admits a global weak solution $u\in L^{\infty}(0,\infty; H_0^1(\Omega))$
with $u_t\in L^2(0,\infty; L^2(\Omega))$ and $u(t)\in W$ for $0\leq t<\infty$.

We first derive the decay rate of $\|u(t)\|_2$. Since $u(t)\in W$ for $0\leq t<\infty$, we have
\begin{equation}\label{2.6}
I(u(t))\geq0,\qquad\ t\geq0.
\end{equation}
By the first equality in \eqref{1.7}, the energy identity \eqref{ei} and \eqref{2.6} we have
\begin{equation}\label{2.7}
\begin{split}
J(u_0)\geq J(u(t))&=\dfrac{a(q-1)}{2(q+1)}\|\nabla u(t)\|_2^2+\dfrac{b(q-3)}{4(q+1)}\|\nabla u(t)\|_2^4+\dfrac{1}{q+1}I(u(t))\\
&\geq \dfrac{a(q-1)}{2(q+1)}\|\nabla u(t)\|_2^2,
\end{split}
\end{equation}
which, together with \eqref{s}, implies that
\begin{equation}\label{2.8}
\|u(t)\|_{q+1}\leq S\|\nabla u(t)\|_2\leq S\Big[\dfrac{2(q+1)}{a(q-1)}J(u_0)\Big]^{\frac{1}{2}}.
\end{equation}
With the help of \eqref{d0}, \eqref{s} and \eqref{2.8} we arrive at
\begin{equation}\label{2.9}
\begin{split}
\|u(t)\|_{q+1}^{q+1}&\leq S^2\|u(t)\|_{q+1}^{q-1}\|\nabla u(t)\|_2^2\\
&\leq S^{q+1}\Big[\dfrac{2(q+1)}{a(q-1)}J(u_0)\Big]^{\frac{q-1}{2}}\|\nabla u(t)\|_2^2\\
&=\Big[\dfrac{J(u_0)}{d_0}\Big]^{\frac{q-1}{2}}a\|\nabla u(t)\|_2^2.
\end{split}
\end{equation}

Taking $\phi=u$ in \eqref{weak-solution} one gets
\begin{equation}\label{2.10}
\dfrac{d}{dt}\|u(t)\|_2^2=-2\Big(a\|\nabla u(t)\|_2^2+b\|\nabla u(t)\|_2^4-\|u(t)\|^{q+1}_{q+1}\Big)=-2I(u(t)).
\end{equation}
Combining \eqref{2.9} with \eqref{2.10} we have
\begin{equation*}
\begin{split}
\dfrac{d}{dt}\|u(t)\|_2^2&\leq -2\Big[a\|\nabla u(t)\|_2^2-\|u(t)\|^{q+1}_{q+1}\Big]\\
&\leq -2a\Big[1-\big(\dfrac{J(u_0)}{d_0}\big)^{\frac{q-1}{2}}\Big]\|\nabla u(t)\|_2^2\\
&\leq -2a\lambda_1\Big[1-\big(\dfrac{J(u_0)}{d_0}\big)^{\frac{q-1}{2}}\Big]\|u(t)\|_2^2,
\end{split}
\end{equation*}
which implies
$$\|u(t)\|_2^2\leq \|u_0\|_2^2e^{-C_1t},$$
where $C_1=2a\lambda_1\Big[1-\big(\dfrac{J(u_0)}{d_0}\big)^{\frac{q-1}{2}}\Big]>0$
and $\lambda_1>0$ is the first eigenvalue of $-\Delta$ in $\Omega$ under homogeneous Dirichlet boundary condition.

Next, we estimate the decay rate of $\|\nabla u(t)\|_2$, $\|u(t)\|_{q+1}$ and $J(u(t))$.
By the definition of $I(u)$ and \eqref{2.9} we have
\begin{equation}\label{2.11}
I(u)\geq a\Big[1-\big(\dfrac{J(u_0)}{d_0}\big)^{\frac{q-1}{2}}\Big]\|\nabla u\|_2^2.
\end{equation}
Define an auxiliary function $H(t)$ by
\begin{equation}\label{2.12}
H(t)=J(u(t))+\|u(t)\|_2^2.
\end{equation}
Then from \eqref{2.7} it follows that
\begin{equation}\label{2.13}
J(u(t))\leq H(t)\leq J(u(t))+\dfrac{1}{\lambda_1}\|\nabla u(t)\|_2^2\leq \Big[1+\dfrac{2(q+1)}{a\lambda_1(q-1)}\Big]J(u(t)).
\end{equation}
Furthermore, the second equality in \eqref{1.7}, \eqref{ei}, \eqref{2.10} and \eqref{2.11} guarantee, for any $\alpha>0$, that
\begin{equation}\label{2.14}
\begin{split}
H'(t)&=-\|u_t(t)\|_2^2-2I(u(t))\\
&\leq-2I(u(t))-\alpha J(u(t))+\dfrac{a\alpha}{4}\|\nabla u(t)\|_2^2+\dfrac{\alpha(q-3)}{4(q+1)}\|u(t)\|_{q+1}^{q+1}+\dfrac{\alpha}{4}I(u(t))\\
&\leq -\alpha J(u(t))+\beta I(u(t)),
\end{split}
\end{equation}
where
$$\beta=-2+\dfrac{\alpha}{4}+\dfrac{\alpha}{4}\Big[1-\big(\dfrac{J(u_0)}{d_0}\big)^{\frac{q-1}{2}}\Big]^{-1}\Big\{1+\dfrac{q-3}{q+1}\big(\dfrac{J(u_0)}{d_0}\big)^{\frac{q-1}{2}}\Big\}.$$
Set $\alpha=8\Big[1-\big(\dfrac{J(u_0)}{d_0}\big)^{\frac{q-1}{2}}\Big]\Big[2-\dfrac{4}{q+1}\big(\dfrac{J(u_0)}{d_0}\big)^{\frac{q-1}{2}}\Big]^{-1}$,
then $\beta=0$ and from \eqref{2.13} and \eqref{2.14} we further obtain
\begin{equation}\label{2.15}
\begin{split}
 H'(t)\leq -\alpha J(u(t)\leq-\alpha\Big[1+\dfrac{2(q+1)}{a\lambda_1(q-1)}\Big]^{-1}H(t):=-C_2H(t).
\end{split}
\end{equation}
where $C_2=\dfrac{a\lambda_1\alpha(q-1)}{a\lambda_1(q-1)+2(q+1)}>0$. Integrating \eqref{2.15} over $[0,t]$ to obtain
\begin{equation}\label{2.16}
J(u(t))+\|u(t)\|_2^2=H(t)\leq H(0)e^{-C_2t}=[J(u_0)+\|u_0\|_2^2]e^{-C_2t}.
\end{equation}
By \eqref{2.7} and \eqref{2.16} we have
\begin{equation}\label{2.17}
\begin{split}
\|\nabla u(t)\|_2^2&\leq \dfrac{2(q+1)}{a(q-1)}J(u(t))\leq \dfrac{2(q+1)}{a(q-1)}[J(u(t))+\|u(t)\|_2^2]\\
&\leq \dfrac{2(q+1)}{a(q-1)}[J(u_0)+\|u_0\|_2^2]e^{-C_2t}.
\end{split}
\end{equation}
Finally, it follows from \eqref{2.9} and \eqref{2.17} that
\begin{equation}\label{2.18}
\|u(t)\|_{q+1}^2\leq S^2\|\nabla u(t)\|_2^2\leq\dfrac{2S^2(q+1)}{a(q-1)}[J(u_0)+\|u_0\|_2^2]e^{-C_2t}.
\end{equation}
The proof of Theorem \ref{rate} is complete. \qquad\qquad\qquad\qquad\qquad\qquad\qquad\qquad\qquad\qquad\qquad\qquad$\Box$

\textbf{Proof of Theorem \ref{bound}.}
We need the following Gagliardo-Nirenberg inequality (see \cite{Brezis}) to derive the lower bound of $\lambda_s$:
\begin{equation}\label{gn}
\|u\|_{q+1}^{q+1}\leq G\|\nabla u\|^{n(q-1)/2}_2\|u\|_2^\gamma,\quad\forall\ u\in H_0^1(\Omega),
\end{equation}
where $G>0$ is a constant depending only on $n$ and $q$ and $\gamma=q+1-\dfrac{n(q-1)}{2}>0$ since $q<2^*-1$.

For any $u\in \mathcal{N}$, by \eqref{gn} we have
\begin{equation}\label{2.20}
a\|\nabla u\|_2^2\leq \|u\|^{q+1}_{q+1}\leq G\|\nabla u\|^{n(q-1)/2}_2\|u\|_2^\gamma,
\end{equation}
which implies
\begin{equation}\label{2.21}
\|u\|_2\geq \Big(\dfrac{a}{G}\Big)^{1/\gamma}\|\nabla u\|_2^{\frac{2}{\gamma}-\frac{n(q-1)}{2\gamma}}
=\Big(\dfrac{a}{G}\Big)^{1/\gamma}\|\nabla u\|_2^{\frac{4-n(q-1)}{2\gamma}},\quad\forall\ u\in \mathcal{N}.
\end{equation}
By the definition of $d$, \eqref{1.7}, \eqref{s} and \eqref{2.20}, it is seen, for any $u\in\mathcal{N}$, that
\begin{equation*}
\begin{split}
d\leq J(u)&=\dfrac{a}{4}\|\nabla u\|_2^2+\dfrac{q-3}{4(q+1)}\|u\|_{q+1}^{q+1}+\dfrac{1}{4}I(u)\\
 &=\dfrac{a}{4}\|\nabla u\|_2^2+\dfrac{q-3}{4(q+1)}\|u\|_{q+1}^{q+1}\\
 &\leq\Big[\dfrac{1}{4}+\dfrac{q-3}{4(q+1)}\Big]\|u\|_{q+1}^{q+1}\leq \dfrac{q-1}{2(q+1)}S^{q+1}\|\nabla u\|_2^{q+1},
 \end{split}
\end{equation*}
which guarantees that
\begin{equation}\label{2.22}
\|\nabla u\|_2\geq\theta:= \Big(\dfrac{2(q+1)d}{q-1}\Big)^{\frac{1}{q+1}}S^{-1},\qquad\forall\ u\in\mathcal{N}.
\end{equation}

If $4-n(q-1)\geq0$, then by combining \eqref{2.21} and \eqref{2.22} we obtain
\begin{equation}\label{2.23}
\lambda_s=\inf\limits_{u\in\mathcal{N}_s}\|u\|_2\geq \inf\limits_{u\in\mathcal{N}}\|u\|_2\geq\Big(\dfrac{a}{G}\Big)^{1/\gamma}\theta^{\frac{4-n(q-1)}{2\gamma}}.
\end{equation}
If $4-n(q-1)<0$, then from \eqref{2.21} and \eqref{52.1} it follows that
\begin{equation}\label{2.24}
\begin{split}
\lambda_s=\inf\limits_{u\in\mathcal{N}_s}\|u\|_2&\geq \Big(\dfrac{a}{G}\Big)^{1/\gamma}\Big[\sup\limits_{u\in\mathcal{N}_s}\|\nabla u\|_2\Big]^{\frac{4-n(q-1)}{2(q+1)-n(q-1)}}\\
&\geq\Big(\dfrac{a}{G}\Big)^{1/\gamma}\Big[\dfrac{2(q+1)s}{a(q-1)}\Big]^{\frac{4-n(q-1)}{4\gamma}}.
\end{split}
\end{equation}

Recalling \eqref{52.1} and Poincar\'{e}'s inequality $\|u\|_2\leq\dfrac{1}{\sqrt{\lambda_1}}\|\nabla u\|_2$ for $u\in H_0^1(\Omega)$
we have
\begin{equation}\label{2.25}
\Lambda_s=\sup\limits_{u\in\mathcal{N}_s}\|u\|_2\leq \sqrt{\dfrac{2(q+1)s}{a\lambda_1(q-1)}}.
\end{equation}
The proof of Theorem \ref{bound} is complete. \qquad\qquad\qquad\qquad\qquad\qquad\qquad\qquad\qquad\qquad\qquad\qquad$\Box$

\textbf{Proof of Theorem \ref{asy}.} Let $u=u(t)$ be a global solution to IBVP \eqref{1.1}.
Without loss of generality, we may assume that
\begin{equation}\label{2.26}
0\leq J(u(t))\leq J(u_0),\qquad t\in[0,\infty).
\end{equation}
In fact, the second inequality follows from \eqref{ei}. If $J(u(t_0))<0$ for some $t_0>0$,
then $I(u(t_0))<0$ by \eqref{1.7}. By Theorem 3.2 in \cite{Han-Li}, $u=u(t)$ blows up in finite time,
which is a contradiction.

Since $J(u(t))$ is non-increasing in $t$ and bounded from below, there exists a constant $J_0\geq0$ such that
\begin{equation}\label{2.27}
\lim\limits_{t\rightarrow\infty}J(u(t))=J_0.
\end{equation}
Letting $t\rightarrow\infty$ in \eqref{ei} and noticing \eqref{2.27} we obtain
\begin{equation*}
\int_0^\infty\|u_{\tau}\|_2^2\mathrm{d}\tau=J(u_0)-J_0\leq J(u_0),
\end{equation*}
which implies that there is an increasing sequence $\{t_k\}_{k=1}^\infty$ with $t_k\rightarrow\infty$ as $k\rightarrow\infty$ such that
\begin{equation}\label{2.28}
\lim\limits_{k\rightarrow\infty}\|u_t(t_k)\|_2=0.
\end{equation}

Denote $u_k=u(t_k)$. We shall show that $\{u_k\}_{k=1}^\infty$ is bounded in $H_0^1(\Omega)$.
For any $v\in H_0^1(\Omega)$, it follows from \eqref{1.5} that
\begin{equation*}
\langle J'(u_k),v\rangle=a(\nabla u_k,\nabla v)+b\|\nabla u_k\|_2^2(\nabla u_k,\nabla v)-(|u_k|^{q-1}u_k,v)=-(u_t(t_k),v),
\end{equation*}
which implies
\begin{equation}\label{2.29}
\begin{split}
\|J'(u_k)\|_{H^{-1}(\Omega)}&=\sup\limits_{\|\nabla v\|\leq 1}|\langle J'(u_k),v\rangle|\leq \sup\limits_{\|\nabla v\|\leq 1}\|u_t(t_k)\|_2\|v\|_2\\
&\leq\dfrac{1}{\sqrt{\lambda_1}}\sup\limits_{\|\nabla v\|\leq 1}\|u_t(t_k)\|_2\|\nabla v\|_2\\
&=\dfrac{1}{\sqrt{\lambda_1}}\|u_t(t_k)\|_2\rightarrow 0,\quad as \ k\rightarrow\infty.
\end{split}
\end{equation}
Therefore, there exists a positive constant $\kappa$ such that
\begin{equation*}
|I(u_k)|=|\langle J'(u_k),u_k\rangle|\leq\|J'(u_k)\|_{H^{-1}(\Omega)}\|\nabla u_k\|_2\leq\kappa\|\nabla u_k\|_2.
\end{equation*}
Recalling \eqref{1.7} again, one gets
\begin{equation*}
\begin{split}
J(u_0)+\dfrac{\kappa}{q+1}\|\nabla u_k\|_2&\geq J(u_k)-\dfrac{1}{q+1}I(u_k)\\
&=\dfrac{a(q-1)}{2(q+1)}\|\nabla u_k\|_2^2+\dfrac{b(q-3)}{4(q+1)}\|\nabla u_k\|_2^4\\
&\geq\dfrac{a(q-1)}{2(q+1)}\|\nabla u_k\|_2^2,
\end{split}
\end{equation*}
which implies that there exists a constant $\Theta>0$ such that
\begin{equation}\label{2.30}
\|\nabla u_k\|_2\leq \Theta,\quad k=1,2,\cdots.
\end{equation}
Therefore, there exists a subsequence of $\{u_k\}_{k=1}^\infty$, which we still denote by $\{u_k\}_{k=1}^\infty$,
and a $u^*\in H_0^1(\Omega)$ such that
\begin{equation}\label{2.31}
\begin{split}
u_k\rightharpoonup u^* \ & weakly \ in \ H_0^1(\Omega) \ as\ k \rightarrow\infty,\\
u_k\rightarrow u^* \ & \ strongly \ in \ L^{q+1}(\Omega) \ as\ k \rightarrow\infty.
\end{split}
\end{equation}

For $u\in H_0^1(\Omega)$, set $E(u)=\dfrac{a}{2}\|\nabla u\|_2^2+\dfrac{b}{4}\|\nabla u\|_2^4$.
Then by Lemma 3.1 in \cite{Han-Li}, $E': H_0^1(\Omega)\rightarrow H^{-1}(\Omega)$ is a strong monotone operator, which satisfies
\begin{equation}\label{2.32}
\langle E'(u)-E'(v),u-v\rangle\geq a\|\nabla u-\nabla v\|_2^2, \quad\forall\ u,v\in H_0^1(\Omega).
\end{equation}
Since $J(u)=E(u)+|u|^{q-1}u$, so
\begin{equation}\label{2.33}
\begin{split}
\langle J'(u_k)-J'(u^*),u_k-u^*\rangle&=\langle E'(u_k)-E'(u^*),u_k-u^*\rangle+(|u_k|^{q-1}u_k-|u^*|^{q-1}u^*,u_k-u^*)\\
&\geq a\|\nabla u_k-\nabla u^*\|_2^2+(|u_k|^{q-1}u_k-|u^*|^{q-1}u^*,u_k-u^*).
\end{split}
\end{equation}
By \eqref{2.29}\ and \eqref{2.30},
\begin{equation}\label{2.34}
|\langle J'(u_k),u_k-u^*\rangle|\leq 2\Theta \|J'(u_k)\|_{H^{-1}(\Omega)}\rightarrow 0,\quad\ as\ k\rightarrow\infty.
\end{equation}
By \eqref{2.31},
\begin{equation}\label{2.35}
|\langle J'(u^*),u_k-u^*\rangle|\rightarrow 0,\quad\ as\ k\rightarrow\infty.
\end{equation}
By H\"{o}lder's inequality, \eqref{2.30} and \eqref{2.31},
\begin{equation}\label{2.36}
\begin{split}
|(|u_k|^{q-1}u_k-|u^*|^{q-1}u^*,u_k-u^*)|&\leq(\|u_k\|_{q+1}^q+\|u^*\|_{q+1}^q)\|u_k-u^*\|_{q+1}\\
&\leq(S^q\|\nabla u_k\|_2^q+\|u^*\|_{q+1}^q)\|u_k-u^*\|_{q+1}\\
&\leq(S^q\Theta^q+\|u^*\|_{q+1}^q)\|u_k-u^*\|_{q+1}
\rightarrow 0,\quad\ as\ k\rightarrow\infty.
\end{split}
\end{equation}
Substituting \eqref{2.34}-\eqref{2.36} into \eqref{2.33} we see that
$$\|u_k-u^*\|_{H_0^1(\Omega)}=\|\nabla u_k-\nabla u^*\|_2\rightarrow0 \quad\ as\ k\rightarrow\infty.$$
Therefore
$$J'(u^*)=\lim\limits_{k\rightarrow\infty}J'(u_k),\qquad in \ H^{-1}(\Omega),$$
which, together with \eqref{2.29}, guarantees that $J'(u^*)=0$, i.e., $u^*\in\mathcal{S}$.
The proof of Theorem \ref{asy} is complete. \qquad\qquad\qquad\qquad\qquad\qquad\qquad\qquad\qquad\qquad\qquad\qquad\qquad\qquad\qquad\qquad\qquad\quad$\Box$

{\bf Acknowledgements}\\
The author would like to express his sincere gratitude to Professor Wenjie Gao for his enthusiastic
guidance and constant encouragement.

\end{document}